Academician Viktor I. Korzyuk[1,2], Jan V. Rudzko[1]

[1]*Institute of Mathematics of the National Academy of Sciences of Belarus, Minsk, Belarus*
[2]*Belarusian State University, Minsk, Belarus*

# CLASSICAL SOLUTION OF THE INITIAL-VALUE PROBLEM FOR A QUASILINEAR WAVE EQUATION WITH DISCONTINUOUS INITIAL CONDITIONS

**Abstract.** For a one-dimensional mildly quasilinear wave equation given in the upper half-plane, we consider the Cauchy problem. The initial conditions have discontinuity of the first kind at one point. We construct the solution using the method of characteristics in an implicit analytical form as a solution of some integro-differential equations. The solvability of these equations, as well the smoothness of their solutions, is studied. For the problem in question, we prove the uniqueness of the solution, and establish the conditions under which its classical solution exists.

**Introduction.** Partial differential equations with discontinuous initial conditions are quite common in various applications, e. g. the propagation of shock waves in a medium [1]. We usually model this phenomenon using the Cauchy problem with discontinuous conditions. This leads to difficulties in the definitions and interpretations of solutions [2].

We often solve such problems by functional methods, e. g. the Fourier transform and the Laplace transform. However, these methods usually do not cover all possible cases of giving the Cauchy conditions [3] since inverse integral transformations converge, as a rule, to the average of the left and right limits [4]. Therefore, various methods have been developed to solve such problems, e. g. the contour integral method [5]. Although this method has a lot of disadvantages, some of which are also inherent in the Fourier method [6], such as increased requirements for smoothness of functions and matching of functions, i.e., functions must satisfy some additional matching and smoothness conditions, it has allowed us to consider a large number of problems of dynamic impact theory, see [7–10] and cited literature. Also, we note the classical d'Alembert method (the method of characteristics), which is not as powerful as functional methods. But we can use it to obtain a qualitative description of impact phenomena [11; 12], to solve some boundary-value problems [13; 14] of impact theory, and it allows us not to identify functions that differ on a set of Lebesgue measure zero.

This work is a continuation of our studies of the Cauchy problem for a mildly quasilinear wave equation [15] and mixed problems with discontinuous initial and boundary conditions [13; 14; 16; 17]. In this article, we consider the Cauchy problem a one-dimensional mildly quasilinear wave equation given in the upper half-plane. The initial conditions of this problem have discontinuity of the first kind at one point of the real axis. We use the method of characteristics to solve this problem. We build the piecewise-smooth solution, which satisfies additional matching conditions, in an implicit analytical form as a solution of some integro-differential equations. We study the solvability of these equations and the smoothness of their solutions. We prove the existence and the uniqueness of the solution to the Cauchy problem under some smoothness conditions of the initial data. The article proposes an approach to constructing solutions with discontinuous initial conditions.

**Statement of the problem.** In the domain $Q = (0, \infty) \times \mathbb{R}$ of two independent variables $(t, x) \in \overline{Q} \subset \mathbb{R}^2$, consider the one-dimensional nonlinear equation

$$\Box u(t,x) + f(t,x,u(t,x),\partial_t u(t,x),\partial_x u(t,x)) = F(t,x),\ (t,x) \in Q, \qquad (1)$$

where $\Box = \partial_t^2 - a^2 \partial_x^2$ is the d'Alembert operator ($a > 0$ for definiteness), $F$ is a function given on the set $\bar{Q}$, and $f$ is a function given on the set $\bar{Q} \times \mathbb{R}^3$. Equation (1) is equipped with the initial condition

$$u(0,x) = \varphi(x), \ \partial_t u(0,x) = \psi(x), \ x \in \mathbb{R}. \qquad (2)$$

where φ and ψ are some real-valued functions defined on the real axis.

The functions φ and ψ are piecewise smooth and defined by the formula

$$\varphi(x) = \begin{cases} \varphi_1(x), & x \in (-\infty, x_0), \\ A, & x = x_0, \\ \varphi_2(x) & x \in (x_0, +\infty), \end{cases} \quad \psi(x) = \begin{cases} \psi_1(x), & x \in (-\infty, x_0), \\ \psi_2(x), & x \in (x_0, +\infty), \end{cases}$$

where $x_0$ and $A$ are arbitrary real numbers, $\varphi_1$ and $\psi_1$ are functions given on the set $(-\infty, x_0]$, and $\varphi_2$ and $\psi_2$ are functions given on the set $[x_0, +\infty)$.

In the case of sufficiently smooth data, namely, $\varphi \in C^2(\mathbb{R})$, $\psi \in C^1(\mathbb{R})$, $F \in C^1(\bar{Q})$, and $f \in C^1(\bar{Q} \times \mathbb{R}^3)$, we considered the problem (1), (2) in the article [15].

For the linear wave equation, i.e., $f \equiv 0$, the problem (1), (2) was studied in the works [16] and [13; 14; 17] in the cases of $\psi \in C^1(\mathbb{R})$ and $\varphi \in C^2(\mathbb{R})$ respectively.

Following [18], we investigate two main questions: 1) the exact statement of the initial value problem (1), (2); 2) the smoothness of the solution. The cause for the first question is the non-existence of the classical solution of the problem due to the condition $\varphi \notin C^2(\mathbb{R})$ or $\psi \notin C^1(\mathbb{R})$.

**Constructing the solution of the Cauchy problem.** We divide the domain $Q$ by the characteristics $x - at = x_0$ and $x + at = x_0$ into three subdomains

$$Q_1 = \{(t,x) | \ (t,x) \in Q \wedge x + at < x_0\},$$
$$Q_2 = \{(t,x) | \ (t,x) \in Q \wedge x - at > x_0\},$$
$$Q_3 = \{(t,x) | \ (t,x) \in Q \wedge x + at > x_0 \wedge x - at < x_0\}.$$

On the closure $\bar{Q}$ of the domain $Q$, we define a function $u$ as the one coinciding with the solution $u^{(j)}$ of the partial differential equation (1)

$$u(t,x) = u^{(j)}(t,x), \quad (t,x) \in Q_*^{(j)}, \qquad (3)$$

on the domain $Q_*^{(j)}$, where $Q_*^{(1)} = \overline{Q^{(1)}} \setminus \overline{Q^{(3)}}$, $Q_*^{(2)} = \overline{Q^{(2)}} \setminus \overline{Q^{(3)}}$, and $Q_*^{(3)} = \overline{Q^{(3)}}$.

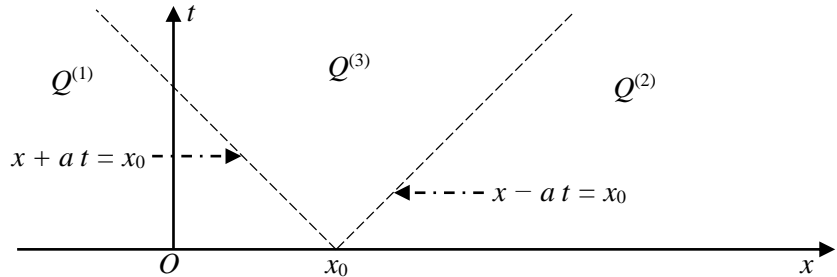

Fig. The partition of the domain Q into three subdomains.

By virtue of discontinuous initial conditions, the problem (1), (2) does not have a classical solution defined on the set $\bar{Q}$. Nevertheless, we can define a classical solution on a smaller set $\bar{Q} \setminus \Gamma$ so that it belongs to the class $C^2(\bar{Q} \setminus \Gamma)$ and satisfies the equation (1), the initial conditions (2), and additional matching conditions given on the set $\Gamma$.

**Definition 1.** A function $u$ is called a classical solution of the problem (1), (2) if the following conditions are satisfied: 1) the function $u^{(j)}$ belongs to the class $C^2\left(\overline{Q_*^{(j)}}\right)$ for each $j \in \{1,2,3\}$; 2) the function $u^{(j)}$ satisfies Eq. (1) in the domain $Q_*^{(j)}$ for each $j \in \{1,2,3\}$; 3) the initial condition $u(0,x) = \varphi(x)$ is met on the whole real axis; 4) the initial condition $\partial_t u(0,x) = \psi(x)$ is met on the whole real axis, except one point $x = x_0$; 5) the Goursat conditions

$$u^{(3)}(t, x = x_0 - at) = u^{(1)}(t, x_0 - at) + A - \varphi_1(x_0), \quad t \geq 0,$$
$$u^{(3)}(t, x = x_0 + at) = u^{(2)}(t, x_0 + at) + A - \varphi_2(x_0), \quad t \geq 0,$$

(4)

hold.

It turns out that finding a classical solution of the problem (1), (2) in the sense of Definition 1 is equivalent to solving the following problem.

**Problem (1), (2) with matching conditions on characteristics.** Find a classical solution of Eq. (1) with the Cauchy conditions (2) and the matching conditions

$$[(u)^+ - (u)^-](t, x = x_0 - at) = A - \varphi_1(x_0), \quad [(u)^+ - (u)^-](t, x = x_0 + at) = \varphi_2(x_0) - A, \quad t \geq 0.$$

Here by $()^\pm$ we have denoted the limit values of the function $u$ and its partial derivatives calculated on different sides of the characteristics $x \pm at = x_0$; i.e., $(\partial_t^p u)^\pm(t, x = r(t)) = \lim_{\delta \to 0+} \partial_t^p u(t, r(t) \pm \delta)$.

The functions $u^{(1)}$ and $u^{(2)}$ are determined from the Cauchy problems

$$\begin{cases} \Box u^{(j)}(t,x) + f(t,x,u^{(j)}(t,x),\partial_t u^{(j)}(t,x),\partial_x u^{(j)}(t,x)) = F(t,x), \ (t,x) \in Q^{(j)}, \\ u^{(j)}(0,x) = \varphi_j(x), \ \partial_t u^{(j)}(0,x) = \psi_j(x), \ (-1)^j x \geq (-1)^j x_0, \end{cases}$$

(5)

for each $j = 1,2$, and under some conditions have the representations

$$u^{(j)}(t,x) = \frac{\varphi_j(x-at) + \varphi_j(x+at)}{2} + \frac{1}{2a}\int_{x-at}^{x+at} \psi_j(\xi)d\xi +$$

$$+ \frac{1}{2a}\int_0^t d\tau \int_{x-a(t-\tau)}^{x+a(t-\tau)} \left(F(\tau,\xi) - f\left(\tau,\xi,u^{(j)}(\tau,\xi),\partial_t u^{(j)}(\tau,\xi),\partial_x u^{(j)}(\tau,\xi)\right)\right) d\xi, \ (t,x) \in \overline{Q_*^{(j)}}.$$

(6)

**Lemma 1.** Let the conditions $F \in C^1(\overline{Q})$, $f \in C^1(\overline{Q} \times \mathbb{R}^3))$, $\varphi_j \in C^2(\mathcal{D}(\varphi_j))$, and $\psi \in C^1(\mathcal{D}(\psi_j))$ be satisfied, and let the function $f$ satisfy the Lipschitz condition with constant $L$ with respect to the three last variables, i. e., $|f(t,x,z_1,z_2,z_3) - f(t,x,w_1,w_2,w_3)| \leq L(|z_1-w_1| + |z_2-w_2| + |z_3-w_3|)$. Then for each $j = 1,2$, the Cauchy problem (5) has a unique solution in the class $C^2(\overline{Q_*^{(j)}})$ determined by the formula (6).

**Proof.** See [15].

Now the function $u^{(3)}$ is determined from the Goursat problem

$$\begin{cases} \Box u^{(3)}(t,x) + f(t,x,u^{(3)}(t,x),\partial_t u^{(3)}(t,x),\partial_x u^{(3)}(t,x)) = F(t,x), \ (t,x) \in Q^{(3)}, \\ u^{(3)}(t, x = x_0 - at) = u^{(1)}(t, x_0 - at) + A - \varphi_1(x_0), \ t \geq 0, \\ u^{(3)}(t, x = x_0 + at) = u^{(2)}(t, x_0 + at) + A - \varphi_2(x_0), \ t \geq 0. \end{cases}$$

(7)

**Lemma 2.** Let the conditions $F \in C^1(\overline{Q})$, $f \in C^1(\overline{Q} \times \mathbb{R}^3))$, $\varphi_j \in C^2(\mathcal{D}(\varphi_j))$, and $\psi \in C^1(\mathcal{D}(\psi_j))$ be satisfied, and let the function $f$ satisfy the Lipschitz condition with constant $L$ with respect to the three last variables, i. e., $|f(t,x,z_1,z_2,z_3) - f(t,x,w_1,w_2,w_3)| \leq L(|z_1-w_1| + |z_2-w_2| + |z_3-w_3|)$. Then the Goursat problem (7) has a unique solution in the class $C^2(\overline{Q_*^{(3)}})$.

**Proof.** According to Lemma 1, the conditions $F \in C^1(\bar{Q})$, $f \in C^1(\bar{Q} \times \mathbb{R}^3))$, $\varphi_j \in C^2(\mathcal{D}(\varphi_j))$, and $\psi \in C^1(\mathcal{D}(\psi_j))$ imply existence and uniqueness of the twice continuously differentiable functions $u^{(1)}$ and $u^{(2)}$. This means that mappings $\gamma_1 : [0,\infty) \ni t \mapsto u^{(1)}(t, x_0 - at) + A - \varphi_1(x_0) \in \mathbb{R}$ and $\gamma_2 : [0,\infty) \ni t \mapsto u^{(2)}(t, x_0 + at) + A - \varphi_2(x_0)$ are twice continuously differentiable too and coincide at the point $t = 0$, i.e., $\gamma_1(0) = \gamma_2(0)$. Now, we use the results of the work [19] to finally prove this lemma.

Let us derive an integro-differential equation for the function $u^{(3)}$. We select four points $A(0, x_0)$, $B\left(\dfrac{x_0 + at - x}{2a}, \dfrac{x_0 - at + x}{2}\right)$, $C(t,x)$, $D\left(\dfrac{at + x - x_0}{2a}, \dfrac{at + x + x_0}{2}\right)$ from the domain $Q_*^{(3)}$, apply the curvilinear parallelogram identity [20] and obtain

$$u^{(3)}(t,x) = u^{(1)}\left(\frac{x_0 + at - x}{2a}, \frac{x_0 - at + x}{2}\right) + A - \varphi_1(x_0) - \varphi_2(x_0) + u^{(2)}\left(\frac{at + x - x_0}{2a}, \frac{at + x + x_0}{2}\right) -$$

$$-\frac{1}{4a^2} \int\limits_{x_*}^{x-at} dy \int\limits_{x_*}^{x+at} \left( F\left(\frac{z-y}{2a}, \frac{z+y}{2a}\right) - f\left(\frac{z-y}{2a}, \frac{z+y}{2a}, u^{(3)}\left(\frac{z-y}{2a}, \frac{z+y}{2a}\right), \right.\right.$$

$$\left.\left. \partial_t u^{(3)}\left(\frac{z-y}{2a}, \frac{z+y}{2a}\right), \partial_x u^{(3)}\left(\frac{z-y}{2a}, \frac{z+y}{2a}\right) \right) \right) dz, \ (t,x) \in Q_*^{(3)}. \quad (8)$$

We state the result as the following assertion.

**Theorem 1.** Let the conditions $F \in C^1(\bar{Q})$, $f \in C^1(\bar{Q} \times \mathbb{R}^3))$, $\varphi_j \in C^2(\mathcal{D}(\varphi_j))$, and $\psi \in C^1(\mathcal{D}(\psi_j))$ hold, and let the function $f$ satisfy the Lipschitz condition with constant $L$ with respect to the three last variables, i. e., $|f(t,x,z_1,z_2,z_3) - f(t,x,w_1,w_2,w_3)| \leq L(|z_1 - w_1| + |z_2 - w_2| + |z_3 - w_3|)$. The initial-value problem (1), (2) has a unique classical solution in the sense of Definition 1. This solution is determined by formulas (3), (6), and (8).

**Analysis of the solution of the Cauchy problem.** Taking into account twice continuous differentiability of the functions $u^{(j)}$ for each $j = 1,2,3$, independence of Lebesgue integral on the behavior of the function on a set of measure zero, the expressions (3), (6), and (8), we can formally rewrite (3) in the form

$$u(t,x) = \frac{\varphi(x-at) + \varphi(x+at)}{2} + \frac{1}{2a}\int\limits_{x-at}^{x+at} \psi(\xi)d\xi + \left(A - \frac{\varphi_1(x_0) + \varphi_2(x_0)}{2}\right)\chi_{Q_*^{(3)}}(t,x) +$$

$$+ \frac{1}{2a}\int\limits_0^t d\tau \int\limits_{x-a(t-\tau)}^{x+a(t-\tau)} \left( F(\tau,\xi) - f(\tau,\xi, u(\tau,\xi), \partial_t u(\tau,\xi), \partial_x u(\tau,\xi)) \right) d\xi, \ (t,x) \in \bar{Q}. \quad (9)$$

where $\chi_A$ is an indicator function of a set $A$. We note that if

$$\left(A - \frac{\varphi_1(x_0) + \varphi_2(x_0)}{2}\right)\chi_{Q_*^{(3)}}(t,x) \equiv 0 \quad (10)$$

then the representation (9) is sometimes called the ('generalized') d'Alembert formula [21].

Following [16], we consider three cases of specifying the discontinuous initial conditions under the conditions of Theorem 1.

1. $\varphi(x_0 - 0) = \varphi(x_0 + 0) = A$, i.e., $\varphi \in C(\mathbb{R})$. From the Goursat conditions (4), and the fact $u^{(j)} \in C^2(\overline{Q^{(j)}})$ for each $j \in \{1,2,3\}$ we see that in this case the solution $u$ belongs the class $C(\bar{Q})$ and satisfies the 'generalized' d'Alembert formula.

2. $\varphi(x_0 - 0) \neq \varphi(x_0 + 0)$ and $A = \dfrac{\varphi(x_0 - 0) + \varphi(x_0 + 0)}{2}$. In this case, the solution $u$ is no longer continuous, but the condition (10) holds, and the solution $u$ satisfies the 'generalized' d'Alembert formula too.

3. $\varphi(x_0 - 0) \neq \varphi(x_0 + 0)$ and $A \neq \dfrac{\varphi(x_0 - 0) + \varphi(x_0 + 0)}{2}$. The solution $u$ is discontinuous and does not satisfy the 'generalized' d'Alembert formula.

These results are consistent with the conclusions of the work [16].

We note that the integro-differential equation (9) can be used to define a mild solution of the problem (1), (2).

**Conclusions.** In the present paper, we have obtained the necessary and sufficient conditions under which there exists a unique classical solution (in an extended sense) of the initial value problem for the mildly quasilinear wave equation with discontinuous initial conditions. And we have proposed an approach to constructing solutions with discontinuous initial conditions, even for nonlinear equations.

**Acknowledgements.** The article was financially supported by the Ministry of Science and Higher Education of the Russian Federation in the framework of implementing the program of the Moscow Center for Fundamental and Applied Mathematics by Agreement no. 075-15-2022-284.

**Information about the authors**

**Viktor I. Korzyuk** – Academician, Professor, Dr. Sc. (Physics and Mathematics), Institute of Mathematics of the National Academy of Sciences of Belarus (11, Surganov Str., 220072, Minsk, Republic of Belarus), Belarusian State University (4, Nezavisimosti Ave., 220030, Minsk, Republic of Belarus). E-mail: korzyuk@bsu.by

**Jan V. Rudzko** – M. Sc. (Mathematics and Computer Sciences). Postgraduate student, Institute of Mathematics of the National Academy of Sciences of Belarus (11, Surganov Str., 220072, Minsk, Republic of Belarus). E-mail: janycz@yahoo.com. https://orcid.org/0000-0002-1482-9106